# The Significance of Ethnomathematics Learning: A Cross-Cultural Perspectives Between Indonesian and Thailand Educators

**I Gusti Ayu Putu Arya Wulandari [1], I Putu Ade Andre Payadnya [1], Kadek Rahayu Puspadewi [1], Sompob Saelee [2]**

[1] Mathematics Education Study Program, Faculty of Teacher Training and Education, Universitas Mahasaraswati Denpasar, Jl. Kamboja No.11A, Dangin Puri Kangin, Kec. Denpasar Utara, Kota Denpasar, Bali Indonesia 80233

[2] Mathematics Education Program, Faculty of Science and Technology, Bansomdejchaopraya Rajabhat University, 1061 Soi Itsaraphap, Itsaraphap Rd. Hiranruchi, Thonburi, Bangkok Thailand 10600

**ABSTRACT**

The field of ethnomathematics holds significance in the pursuit of comprehending how students can grasp, express, manipulate, and ultimately apply mathematical concepts. However, ethnomathematics is also considered a complex concept in Asian countries such as Indonesia and Thailand, which can pose challenges as it needs to be comprehensively understood. This research aims to fill the gap by understanding the cross-cultural perspective of mathematics educators in Indonesia and Thailand. The participants were lecturers, teachers, and pre-service teachers. Data was gathered through questionnaires and interviews. The analytical approach involved were data reduction, data presentation, drawing conclusions or verification, and data validity. Positive responses were indicated by mathematics educators with the average scores of respondents in Indonesia at 4.77 and Thailand at 4.57. This research concludes the importance of integrating ethnomathematics in education, which is closely tied to cultural development, emphasizing the crucial role of employing comprehensive strategies in its implementation.

## INTRODUCTION

Apart from its association with the progress of science, mathematics is also linked to cultural advancement. Culture-based mathematics, commonly known as ethnomathematics, has become a trend in modern mathematics education. Ethnomathematics is characterized as the mathematical practices observed within cultural groups, encompassing urban and rural communities, labor groups, children of specific age groups, indigenous communities, and others (Rahmawati, 2012). This trend in education can be viewed as a program designed to explore how students can not only understand but also express, manipulate, and ultimately apply mathematical ideas, concepts, and practices to address problems rooted in their daily activities.

Indonesia and Thailand are countries in Southeast Asia known for their cultural richness. Due to the diversity of cultures, the practice of ethnomathematics has also flourished in these two nations. Prahmana and D'Ambrosio (2020) discovered that in Indonesia, ethnomathematics extends beyond merely exploring and experimenting with cultural elements in mathematics education within various schools. Their findings suggest that ethnomathematics can be integrated into the formal mathematics education curriculum in Indonesia. For instance, Suryanatha and Apsari studied the concept of modulo found in the naming system in Bali (Suryanatha & Apsari, 2022). Furthermore, Puspadewi and Putra examined ethnomathematics in Balinese weaving crafts and its connection to learning (Puspadewi & Putra, 2014). From weaving crafts, individuals can acquire insights into various mathematical concepts, such as principles related to tiling, parallel lines, and angles.

The exploration of ethnomathematics in Thailand has also entered the realm of school education. For example, Suryana et al. (2022) explored ethnomathematics in the five tones of the Thai language. They found that the Thai language tones can be studied as an ethnomathematics object in mathematics learning, especially in the topics of curves and functions. This demonstrates that the exploration and application of culturally based mathematics in education are feasible in both Indonesia and Thailand. However, there are undoubtedly barriers and challenges in practice.

The implementation of ethnomathematics in education has undergone various developments and challenges. This is due to the fact that ethnomathematics is not easily embraced by students. Rosa and Orey expressed concerns about the success of ethnomathematics as a pedagogical action (Rosa & Orey, 2011). This arises from the lack of inclusion of ethnomathematics in classroom materials, insufficient discussions about ethnomathematics, particularly in higher education, inadequate training of teachers to understand the wide application of ethnomathematics, and the shallow understanding of ethnomathematics activities. Additionally, many teachers feel that students will not learn mathematics according to the demands of the traditional curriculum if they employ an ethnomathematical approach. The ethnomathematics perspective is believed to make the accountability of students and teachers in national exams problematic. Furthermore, Nur et al. (2019) added that the concept of ethnomathematics is a broad and complex idea that must be fully understood by every educator. Hence, the teacher plays a pivotal role in the success of ethnomathematics-based learning.

Such contradictions can pose obstacles to the implementation of ethnomathematics in school education, leading to suboptimal ethnomathematics learning that results in students' low academic achievement. Nonetheless, ethnomathematics plays a pivotal role in enhancing students' mathematical abilities. Richardo highlighted the outcomes of his research, revealing that the integration of ethnomathematics into mathematics education introduces a new dimension, suggesting that learning mathematics extends beyond the classroom. Engaging with the outside world, visiting, or interacting with local culture can serve as valuable tools for learning mathematics (Richardo, 2016). Furthermore, in the era of the fourth industrial revolution, ethnomathematics serves as a significant bridge between the advancements in information and communication technology and the preservation of Indonesian culture across generations (Darmayasa et al., 2019).

Given the importance of ethnomathematics in education and the various challenges in its implementation, especially in culturally rich countries like Indonesia and Thailand, researchers argue that it is important to conduct research on the intercultural perspective on the significance of ethnomathematics learning between Indonesian and Thailand educators. Therefore, this research, building upon the most recent scholarship, significantly fills the gap in previous studies on ethnomathematics in the context of Indonesia and Thailand. While previous research, such as that of Prahmana & D'Ambrosio (2020), has acknowledged the presence and potential of ethnomathematics in these countries, a comprehensive comparative study examining the perspectives of educators from both nations regarding the significance and implementation of ethnomathematics is notably absent. By delving into the intercultural perspective on the importance of ethnomathematics learning, this research addresses the critical need to understand the cultural nuances that shape the reception and application of ethnomathematics in educational settings. This study's novelty lies in its explicit focus on the viewpoints of educators, with the intention of offering a comprehensive understanding of the challenges, opportunities, and potential strategies for integrating ethnomathematics into the curriculum, thereby fostering more effective and culturally responsive practices within mathematics education.

In the subsequent section, we delve into the methodology employed in this research, detailing the research location and participants, research design and procedure, data collection techniques, and data analysis techniques. The third section encompasses the outcomes of the questionnaire score analysis, as well as the discussion of the comprehensive interviews conducted with representative respondents, revealing their responses to the significance of ethnomathematics in learning, along with various opinions and inputs from the respondents. Finally, the conclusion outlines the research's findings, limitations, and recommends future studies.

**RESEARCH METHOD**
**Research Location and Participants**

The research was conducted at Mahasaraswati University Denpasar and Bansomdejchaopraya Rajabhat University, Thailand, along with secondary schools affiliated with both institutions. The research took place during the Even Semester of the Academic Year 2022/2023. The participants in this research were mathematics educators consist of: lecturer, high school teachers, secondary school teachers, and pre-service teachers from Indonesia (138 participants) Thailand (145 participants).

**Research Design and Procedure**
1) Research Design
   This study employed a descriptive qualitative research method aiming to provide a more detailed perspective on the intercultural aspects of mathematics educators concerning the significance of ethnomathematics in mathematics education. Qualitative research is intended to understand phenomena such as the experiences of research subjects, including behaviors, perceptions, motivations, actions, and more. This is done holistically, described through words and language, in a specific natural context, and utilizing various scientific methods (Semiawan, 2010). The qualitative approach was chosen to reveal a deeper understanding of the intercultural perspective of mathematics educators regarding the significance of ethnomathematics in mathematics education. The study will employ a qualitative descriptive research approach, gathering data directly from questionnaire responses and interviews.
2) Research Procedure
   The research procedure encompasses the activities designed by the researchers to be implemented during the research. The steps involved in this research procedure are as follows: a) Determining the research steps through a focus group discussion with partner institutions, b) Selecting research participants from both institutions, b) Establishing the analysis procedures for examining the participants' perspectives, d) Preparing research instruments such as questionnaires and interview guidelines, e) Consulting research instruments with the research team and experts, f) Conducting instrument pilot tests to ensure their validity and reliability, g) Implementing the research by distributing instruments and conducting interviews, h) Examining and analyzing the acquired data, and i) Drawing conclusions from the data analysis.

**Data Collection Techniques**

Several data collection methods were utilized in this research, including the following:

1) Questionnaire Technique

    The questionnaire is a data collection method used to understand individuals. It consists of a list of questions on various aspects. According to Sugiyono (2014), a questionnaire is a data collection technique where the researcher provides a list of written questions or statements to be answered by respondents. In this study, a questionnaire was used to obtain intercultural perspectives on the significance of ethnomathematics learning among educators from Indonesia and Thailand. The questionnaire, consisting of 20 questions with score 1-5, specifically designed to assess educators' perspectives, was created using Google Forms and distributed through WhatsApp and Line.

2) Interview Technique

    Interviews involve obtaining information or data for research purposes through a question-and-answer process, conducted face-to-face between the interviewer and the respondent using an interview guide (Siregar, 2013). Interviews were conducted to reconfirm and reinforce the findings obtained in the research. Through interviews researchers can obtain various answers and points of view which are important for the richness of research data. Interviews in this research were conducted with 8 participants consisting of junior high school teacher, high school teacher, mathematics education lecturer, and pre-service teacher in Indonesia and Thailand. Interviews were carried out directly by first selecting representative respondents as representatives by looking at the results of the questionnaire. The type of interview used is an unstructured interview which focuses on digging deeper into the opinions of mathematics educators regarding the significance of culture-based mathematics learning.

3) Documentation

    Documentation involves gathering data on variables from various sources, including records, transcripts, books, student grade lists, attendance records, and more. The documentation method was used to acquire data regarding the implementation of ethnomathematics teaching in the institutions involved in this research.

**Data Analysis Technique**

The research employed a qualitative descriptive data analysis technique, with the following steps:

1) Data Reduction

    Data reduction involves sharpening, categorizing, directing, discarding unnecessary data, and organizing the remaining data to draw and verify final conclusions. This process includes selecting, focusing, simplifying, and abstracting the raw data written in field notes. The data reduction steps in this research were as follows: a) Examining the questionnaire results and grouping them based on participant responses, b) Transforming the questionnaire data obtained from the participants into notes for use in the interviews, c) Simplifying the conducted interview results into a well-structured and organized form, and then transforming them into notes. In grouping the questionnaire results, positive responses were obtained from an average of agree and strongly agree responses with a score of 4-5.

2) Data Presentation

    Data presentation involves organizing a set of structured information that allows for drawing conclusions and taking action. In this stage, the questionnaire results were organized according to the participants.

3) Drawing Conclusions or Verification

    Verification is part of a comprehensive configuration activity that can answer research questions and objectives. Comparing the questionnaire results with the interview outcomes can lead to conclusions about the educators' perspectives.

4) Data Validity Check

    After analyzing the available data to find answers to the research problem, the next step involves examining the validity of the findings. To determine the validity of the findings, a verification technique is required. In this research, data validity checking employed the triangulation technique.

The triangulation technique involves using external sources of data to check or compare against the data. The type of triangulation used in this research was source triangulation, comparing and cross-checking the degrees of belief of information obtained through different times and methods in the qualitative approach. The source triangulation stage conducted in this research involved comparing the results of the questionnaire with the interview outcomes.

**RESULTS AND DISCUSSION**

Dissemination of the questionnaire was carried out among various respondents from the ranks of university lecturers, high school teachers, junior high school teachers, as well as pre-service teachers from various educational institutions in Indonesia and Thailand. From the distributed questionnaire results, a score recapitulation was obtained, showing the responses of educators from Indonesia and Thailand regarding the significance of ethnomathematics in learning.

**Table 1. Recapitulation of Average Scores of Indonesian Respondents**

| No | Respondent | Indonesia | | | Thailand | | |
|---|---|---|---|---|---|---|---|
| | | Number of Respondent | Score | Criteria | Number of Respondent | Score | Criteria |
| 1 | Lecturer | 30 | 4,90 | Positive | 10 | 4,75 | Positive |
| 2 | High School Teacher | 55 | 4,53 | Positive | 50 | 4,62 | Positive |
| 3 | Middle School Teacher | 40 | 4,78 | Positive | 28 | 4,10 | Positive |
| 4 | Pre-Service Teacher | 20 | 4,85 | Positive | 50 | 4,80 | Positive |
| | Total | 145 | Average = 4,77 | | 138 | Average = 4,57 | |

Table 1 presents a comprehensive recapitulation of the average scores from Indonesian respondents, including 30 lecturers, 55 high school teachers, 40 middle school teachers, and 20 pre-service teachers. The data reveals an overall high average score of 4.77, indicating a strong acknowledgment of the significance of ethnomathematics within the Indonesian educational context. The scores signify a positive attitude and recognition of the importance of incorporating cultural elements into mathematics education, as evidenced by the consistently high ratings across all categories of respondents. These findings align with recent research emphasizing the benefits of integrating ethnomathematics in diverse educational settings, fostering a deeper understanding of mathematical concepts through the lens of cultural relevance (Brandt & Chernoff, 2014).

Similarly, Table 1 also illustrates the average scores obtained from Thai respondents, comprising 10 lecturers, 50 high school teachers, 28 middle school teachers, and 50 pre-service teachers. The table reflects a slightly lower overall average score of 4.57, suggesting a comparable but marginally more reserved stance towards the integration of ethnomathematics within the Thai educational landscape. While the scores still indicate a high level of acknowledgment, the subtle differences in the averages might point to varying perspectives on the extent to which cultural elements should be integrated into mathematics education. Recent studies have highlighted the significance of considering cultural context in mathematics instruction, emphasizing the potential for enhancing students' engagement and understanding through culturally relevant pedagogical practices (Acharya, et al, 2021). These observations underscore the importance of further research and collaborative efforts to bridge the gap between cultural diversity and mathematical education, fostering inclusive learning environments that cater to the multifaceted needs of students within diverse sociocultural contexts.

Interviews were conducted with one representative from each respondent group. The results of interviews with Indonesian and Thai educators show similar perceptions regarding the significance of ethnomathematics in learning. The following is an excerpt from an interview with one of the middle school teacher representatives.

**Researcher (R):** Good morning, thank you for taking the time to discuss the role of ethnomathematics in contemporary education with us. Do you think ethnomathematics is important for mathematics learning today?

**Teacher (T):** Good morning, it's my pleasure. I believe ethnomathematics is crucial in today's learning environment as culture continues to evolve.

**R:** Indeed, your viewpoint aligns with recent literature emphasizing the importance of integrating cultural perspectives into education. Could you elaborate on how ethnomathematics can foster a deeper appreciation of cultural diversity among students?

**T:** Certainly, by incorporating diverse cultural perspectives into mathematical instruction, we not only enhance students' understanding of different cultures but also instill a sense of inclusivity within the learning process.

**R:** Your perspective echoes recent research, highlighting the need for a delicate balance in bridging cultural perspectives and rigorous high school mathematics. Moving on, how do you perceive the relationship between ethnomathematics and realistic mathematics learning in middle-school?

**T:** Ethnomathematics is closely intertwined with realistic mathematics learning. It helps students understand mathematical concepts in real-life contexts, fostering a deeper connection between theoretical knowledge and practical applications.

**R:** Your insights emphasize the role of ethnomathematics in promoting authentic learning experiences for junior high school students. Lastly, based on your experience, how do you think ethnomathematics can shape the future of meaningful mathematics education?

**T:** Ethnomathematics has the potential to revolutionize the way we approach mathematics education. By incorporating culturally relevant teaching practices, we can stimulate students' curiosity and engagement, fostering a positive learning environment that celebrates cultural diversity within mathematical concepts.

**R:** Your valuable insights shed light on the transformative role of ethnomathematics in shaping the future of meaningful mathematics education. Thank you for sharing your expertise with us today.

**T:** You're welcome. It was a pleasure discussing this vital topic with you.

The insights provided by the lecturer resource person offer a compelling perspective on the enduring significance of ethnomathematics in contemporary learning environments. According to the interviewee, the dynamic nature of culture necessitates the integration of ethnomathematics as a crucial component of the educational process. The viewpoint underscores the evolving nature of cultural contexts, emphasizing the need for educators to incorporate diverse cultural perspectives into mathematical instruction. This aligns with recent scholarly literature, which emphasizes the importance of culturally responsive pedagogy in enhancing students' engagement and understanding of mathematical concepts within diverse societal frameworks (Anhalt et al., 2018). The lecturer's assertion further underscores the role of ethnomathematics in fostering a deeper appreciation of cultural diversity, paving the way for inclusive educational practices that resonate with students from various cultural backgrounds.

The feedback from the high school teacher resource persons presents a nuanced understanding of the practical challenges associated with integrating ethnomathematics into complex high school curricula. While acknowledging the importance of the concept in enhancing students' problem-solving abilities, the interviewees express concerns about the feasibility of seamlessly incorporating ethnomathematics into the advanced high school material. This sentiment echoes recent research that highlights the intricate balance required to bridge cultural perspectives with complex mathematical concepts within the high school educational context (Nasir et al., 2008). The teachers' feedback suggests the need for comprehensive strategies that facilitate the effective integration of ethnomathematics, ensuring its seamless alignment with the rigorous demands of high school mathematics education.

The perspectives shared middle- school teacher resource persons emphasize the interconnected nature of ethnomathematics and realistic mathematics learning. The interviewee underscores the role of ethnomathematics in facilitating students' comprehension of mathematical concepts within real-world contexts. This sentiment echoes contemporary research highlighting the significance of real-life applications in fostering students' conceptual understanding and problem-solving skills in mathematics (Laurens et al., 2018). The junior high school teachers' feedback underscores the potential of ethnomathematics as a catalyst for promoting authentic learning experiences, nurturing students' ability to perceive the practical relevance of mathematical concepts in their immediate surroundings. This insight accentuates the transformative potential of ethnomathematics in cultivating students' holistic mathematical understanding within the junior high school curriculum.

The perspective shared by the pre-service teacher sheds light on the promising role of ethnomathematics as a catalyst for meaningful mathematics learning. The interviewee endorsement of ethnomathematics as the future of enriching mathematical education underscores the transformative potential of culturally responsive pedagogical practices. This sentiment resonates with recent research emphasizing the role of culturally relevant teaching approaches in fostering students' intrinsic motivation and interest in mathematics (Roche et al., 2021). The student optimistic outlook on ethnomathematics highlights its capacity to stimulate students' curiosity and engagement, fostering a positive learning environment that nurtures their appreciation for the cultural diversity embedded within mathematical concepts. This underscores the need for educators to embrace innovative pedagogical strategies that integrate ethnomathematics, fostering a dynamic learning ecosystem that resonates with students' evolving educational needs and aspirations.

The interview findings also showed that the application of ethnomathematics in the educational context differs between Indonesia and Thailand. While both countries recognize the significant role of ethnomathematics in fostering cultural inclusivity and enhancing students' understanding of mathematical concepts, the interviews with respondents from Thailand suggest a limited integration of ethnomathematics within the learning process. Despite acknowledging its importance, the Thai educational landscape appears to demonstrate a relatively lower prevalence of ethnomathematics in instructional practices. This observation is in line with recent research emphasizing the challenges associated with implementing ethnomathematics in diverse educational settings, particularly in contexts where traditional pedagogical approaches remain prevalent (Jun-on & Suparatulatorn, 2023). In contrast, the interviews with Indonesian respondents reflect a more proactive approach to integrating ethnomathematics, highlighting its relevance and practical application within the Indonesian educational framework (Mania & Alam, 2021). This discrepancy underscores the need for tailored educational policies and collaborative efforts to promote the effective integration of ethnomathematics in diverse cultural contexts, fostering inclusive learning environments that cater to the multifaceted needs of students within the Southeast Asian region.

Based on the comprehensive discussions with various educational stakeholders, it is evident that the integration of ethnomathematics within the educational framework holds promise for fostering a culturally inclusive and engaging learning environment. While emphasizing the importance of embracing diverse cultural perspectives in mathematics education (Gutierrez, 2017), the discussions highlighted the challenges associated with seamlessly integrating ethnomathematics into advanced high school curricula (Freire & McCray, 2020). Moreover, the insights underscored the transformative potential of ethnomathematics in promoting authentic learning experiences, fostering a deeper connection between theoretical knowledge and practical applications within the junior high school curriculum (Rosa & Orey, 2015). By embracing culturally relevant teaching practices, ethnomathematics has the capacity to stimulate students' curiosity, engagement, and appreciation for the cultural diversity embedded within mathematical concepts, envisioning a more inclusive and dynamic educational landscape (Nasir, 2021).

**CONCLUSIONS**

This research concludes positive responses that were demonstrated by mathematics educators regarding the significance of ethnomathematics in education, with an average score of respondents in Indonesia and Thailand reach 4.77 and 4,57 respectively, indicating strong recognition of the importance of ethnomathematics in the Indonesian educational context. The study also concluded that ethnomathematics, in integrating culture into the teaching of mathematics in Indonesia and Thailand, despite challenges in material integration, especially at higher levels, has great potential to enhance students' understanding of mathematical concepts in everyday life. This research can serve as a foundation and consideration in the development of culturally-based ethnomathematics curriculum that is appropriate and applicable across all levels of education.

One of the primary limitations of this research is the constrained number of respondents and participating institutions. The pool of respondents was restricted to educational institutions in Bali, its surrounding regions, and Thailand. Notably, the participating schools primarily consisted of affiliated or collaboratively associated institutions with the research entity. These constraints may impact the generalizability and broader applicability of the study's findings. The subsequent research is expected to involve a larger number of respondents from diverse cultural backgrounds. This is anticipated to enhance the understanding of various cultural contexts associated with mathematics education. By incorporating cross-cultural perspectives, the research will enable the identification of different cultural patterns in the understanding and application of mathematical concepts. Consequently, the study will make a more significant contribution in emphasizing the importance of cultural integration in an inclusive and sustainable approach to mathematics education.

**CONFLICT OF INTEREST**

The authors declare there is no conflict of interest.

**ACKNOWLEDGEMENTS**

The authors would like to thank Universitas Mahasaraswati Denpasar, Bali, Indonesia for funding this research through an Overseas Collaboration Research Grant through research contract number: K.127/B.01.01/LPPM-Unmas/IV/2023. The author also would like to thank Bansomdejchaopraya Rajabhat University, Thailand for the excellent collaboration in carrying out this research.